\documentclass[letterpaper, 10 pt, conference]{ieeeconf}  

\IEEEoverridecommandlockouts                              

\usepackage{graphics} 
\usepackage{epsfig} 
\usepackage{mathptmx} 
\usepackage{times} 
\usepackage{amsmath} 
\usepackage{amssymb}  
\usepackage{verbatim}
\usepackage{subfig}
\usepackage{setspace}
\usepackage{algorithm}
\usepackage{algpseudocode}
\usepackage{float}
\title{\LARGE \bf
Bayesian Persuasive Driving
}

\author{Cheng Peng and Masayoshi Tomizuka
\thanks{This work was supported by China Scholarship Council (CSC) Scholarship.}
\thanks{C. Peng and M. Tomizuka are with the Department of Mechanical
Engineering, University of California, Berkeley, CA 94720 USA (e-mail:
chengpeng2014@berkeley.edu; tomizuka@berkeley.edu).}%
}

\begin{document}

\maketitle
\thispagestyle{empty}
\pagestyle{empty}

\begin{abstract}
In the autonomous driving area, interaction between vehicles is still a piece of puzzle which has not been fully resolved. The ability to intelligently and safely interact with other vehicles can not only improve self driving quality but also be beneficial to the global driving environment. In this paper, a Bayesian persuasive driving algorithm based on optimization is proposed, where the ego vehicle is the persuader (information sender) and the surrounding vehicle is the persuadee (information receiver). In the persuasion process, the ego vehicle aims at changing the surrounding vehicle's posterior belief of the world state by providing certain information via signaling in order to achieve a lower cost for both players. The information received by the surrounding vehicle and its belief of the world state are described by Gaussian distributions. Simulation results in several common traffic scenarios are provided to demonstrate the proposed algorithm's capability of handling interaction situations involving surrounding vehicles with different driving characteristics.
\end{abstract}


%
\IEEEpeerreviewmaketitle

\section{Introduction}
Although autonomous vehicles have been spotted more and more frequently driving on the city roads, most of them are still not interacting with other road users like human drivers do. Instead, most autonomous vehicles are implementing a reactive behavior, which means that the trajectory predictions about surrounding vehicles are made first and the ego vehicle's driving actions are decided accordingly by applying obstacle avoidance algorithm. However, in this planning pattern, the interaction between vehicles is ignored since the interacting vehicle's future driving profile is assumed to be independent of the ego vehicle's behavior. Therefore, more efforts are needed to fulfill a real interactive, efficient and cooperative driving environment where robot cars and human drivers coexist.

Several approaches have been proposed for interactive driving in literatures, which can be categorized into two groups in general. The first category is multi-agent algorithm \cite{de2013autonomous}-\cite{hafner2013cooperative}, where the assumption is made that all the vehicles involved in a driving scenario can be controlled. In a multi-agent system, all the vehicles optimize the same functional and each vehicle knows exactly what the others will be doing and what influence its behavior will cause. The major drawback of the multi-agent algorithm is that in the real world, not all the related vehicles can be controlled, especially at the current stage where autonomous vehicles only make a minority of all road users. Another disadvantage is its heavy dependency on vehicle-to-vehicle (V2V) communication which may not be reliable enough in real system. 

The second category of algorithm is based on interactive prediction of surrounding vehicle's future behaviors. Several promising prediction approaches have been proposed in literatures including partially observable Markov decision process (POMDP), deep neural network, optimization-based method and so on. For example, a POMDP based decision making strategy was proposed in \cite{hubmann2017decision} for intersection scenario where interacting vehicles were assumed to pick one route from a predetermined route hypothesis set. However, this approach was confined to a specific driving situation due to the fixed route hypothesis setting. In \cite{lenz2017deep}, several deep neural network based motion models were evaluated for the highway entrance scenario, among which a model was selected for fast computation and relatively good performance. The drawback of this approach is that its performance heavily depends on feature selection which is not general for different driving scenarios. Optimization provides another direction of generating interactive prediction for surrounding vehicles. As an example, in \cite{sadigh2016planning}, the human driver behavior is predicted by optimizing a reward function with pre-defined structure and parameters learned via inverse reinforcement learning (IRL) algorithm. Through including the learned human model reward function in robot vehicle trajecotry planning, the interaction effects can be handled. Compared to the POMDP and neural network based methods, this optimization-based approach is more computational efficient. However, the learned reward function can only describe a particular type of driver, which is not a general solution to the interactive driving problem.

Another optimization-based interactive planning approach is to formulate the interactive driving problem as a Bayesian persuasion game, which is first proposed for economic application \cite{kamenica2011bayesian}. In the Bayesian persuasion game, there is one sender with information who attempts to persuade the receiver to change his/her action so that the welfare of both players can be improved. The basic assumptions include 1) the receiver's behavior is dependent on his/her belief of the world state and 2) both the players are rational Bayesian under which the interaction can be described as a Bayesian process. The persuasion process can be achieved by the sender via selecting certain information to convey to the receiver so that the receiver's posterior belief distribution of the world state can be properly manipulated.

In this paper, the world state for the interactive driving environment is defined to be the ego vehicle's conservativeness perceived by the surrounding vehicle, the ego vehicle is defined as the information sender, whose information of driving intention can be reflected from his/her driving behavior and the surrounding interacting vehicle is the corresponding information receiver. With regard to the signal, several candidates are available including binary signal of yield or not yield, discrete signal of driving route selection and continuous signal of driving state. As a starting point, the ego vehicle's continuous driving state is selected as the signal since it carries the most driving information. By determining the optimal signal based on an optimization, the ego vehicle is able to achieve maximization/minimization of utility/cost expectation for both players. Meanwhile, the receiver will extract more information about the world state from the perceived signal and thus his/her posterior belief of the world state can be updated.

The remainder of the paper is structured as follows. In Section \ref{sec:original_bayes}, the mathematical formulation of the general Bayesian persuasion problem is introduced. In Section \ref{sec:bay_per_car}, the concrete Bayesian persuasion problem for interactive driving is formulated as an optimization with certain constraints incorporated including vehicle dynamics, safety and physical saturation limit. In addition, the integrals are calculated and the optimization problem is reformulated into a tractable form. In Section \ref{sec:simu}, the proposed algorithm's effectiveness is illustrated by simulations in several driving scenarios. Section \ref{sec:conc} concludes the paper.

\section{General Bayesian Persuasion}\label{sec:original_bayes}
In this section, the formulation of a general Bayesian persuasion game is introduced and the related preliminary notations are given. 

A general and intuitive definition of the persuasion problem is to exploit some information advantage to influence the opponent's action or intention. In fact, the persuasion behavior is ubiquitous in everyday life with applications in a great deal of areas including economy, psycology, decision making theory and so on. Basically, in almost any interacting process, there always exists a persuasion scheme which is advantageous to some or all players. 

Among various persuasion models, the Bayesian persuasion model \cite{kamenica2011bayesian} first proposed for economy application stands out as the most popular and fundamental one. In a Bayesian persuasion game, there are two players: the player sending information is called sender and the other one receiving information is called receiver. The sender aims to change the receiver's action so that there is a higher probability that a situation more beneficial to both parties can be achieved. The receiver's task is to pick an action based on the information extracted from the sender's signal. The reward/cost function of the game depends on both players' actions meaning that neither player can determine the game's result by himself. This characteristic leads to the fact that both players cannot exactly know what the reward or cost will be until both of their actions have been unveiled.

In order to formulate the persuasion game mathematically, some notations are introduced first. The action of the receiver is denoted as $a\in A$, where $A$ is the receiver's action space. The world state is represented by $\omega\in \Omega$, where $\Omega$ is the world state space. The realization space of the sender's information signal is denoted as $S$ and the corresponding signal realization is $s$. However, the information carried by the signal may not be fully comprehended by the receiver. In the receiver's point of view, the signal should be described by a probability distribution which is called signal belief distribution denoted by $\pi(s)$. Hence, $\pi(s)$ represents the distribution of the information signal perceived by the receiver. Due to the similar reason, another world state belief distribution $\mu(\omega|s)$ is introduced to describe how the world state $\omega$  is influenced by the signal in the receiver's mind. For the game's objective function, maximizing a reward function and minimizing a cost function are equivalent. Thus, without loss of generality, we choose to define a cost function $c(\omega, a, s)$, which is dependent on the world state, receiver's action and the sender's signal. With the notations defined above, the Bayesian persuasion problem can be formulated as:
\begin{eqnarray}\label{eq:bay_per}
\min_{\pi(s)} & E_{\pi(s)}\hat{c}(\mu(\omega|s)),
\end{eqnarray}
where
$$\hat{c}(\mu(\omega|s))=\min_{a}E_{\mu(\omega|s)}c(\omega,a,s)$$
represents the expected cost at a specific signal realization $s$ when the receiver holds the belief distribution $\mu(\omega|s)$. The solution to the optimization (\ref{eq:bay_per}) is defined as the optimal signal and the corresponding achieved minimum value is called the value of the optimal signal. In the Bayesian persuasion process described by (\ref{eq:bay_per}), according to the rationality assumption, the receiver will decide his action by optimizing the objective function expectation given his belief of the world state influenced by the sender's signal. Moreover, since the sender is aware that the receiver's rational, he can then determine the optimal signal to send based on the receiver's strategy.

In summary, the optimization problem in (\ref{eq:bay_per}) illustrates the idea that the sender's purpose is to minimize both players' cost by manipulating the receiver's posterior belief distribution $\mu(\omega|s)$ via conveying information via signaling perceived by the receiver as a probability distribution $\pi(s)$.

\section{Bayesian persuasive driving}\label{sec:bay_per_car}
In this section, the Bayesian persuasion framework is applied to the interacting driving problem with Gaussian distribution assumption. The concrete definition of variables are given in the autonomous driving context first and a mathematical approximation is then applied to the resulting driving persuasion problem to make it tractable.

\subsection{Bayesian persuasion in autonomous driving context}
In order to formulate a Bayesian persuasion game in the autonomous driving background, the players along with their possible actions and signals need to be clearly defined first.

Intuitively, the surrounding interacting vehicle is selected as the information receiver and his driving state
$$\textbf{x}_t^s=[\begin{array}{cccc}x^s_t & y^s_t & \theta^s_t & v^s_t\end{array}]^T$$
is defined as his action $a_t$, where $x^s_t$, $y^s_t$, $\theta^s_t$ and $v^s_t$ denote the surrounding vehicle's $x-y$ positions, yaw angle and speed in the lane-based coordinate frame at time instant $t$ respectively. The information sender role is then naturally assigned to the ego vehicle. With regard to the signaling content, there are quite a few options including intention indicator for yielding or not yielding, route selection preference and the driving state itself. In this paper, the ego vehicle's driving state
$$\textbf{x}_t^e=[\begin{array}{cccc}x^e_t & y^e_t & \theta^e_t & v^e_t\end{array}]^T$$
is chosen as the signal realization $s_t$ since it includes more intention/driving behavior information and it is directly perceivable for the surrounding vehicle. Similar with the notation for the surrounding vehicle, $x^e_t$, $y^e_t$, $\theta^e_t$ and $v^e_t$ denote the ego vehicle's $x-y$ positions, yaw angle and speed in the lane-based coordinate frame at time instant $t$ respectively.

The last definition involved in the Bayesian persuaison is the world state variable $\omega$. The characteristic of the Bayesian persuasion game imposes two requirements on the state $\omega$, which are 1) it determines the cost for both players along with the receiver's action and the sender's signal, 2) it cannot be directly influenced by the receiver's action. In order to satisfy the two mentioned properties, the state is defined as the ego vehicle's conservativeness perceived by the surrounding vehicle. Basically, the state $\omega$ denotes the receiver (surrounding vehicle)'s impression of the sender (ego vehicle), whether aggressive, conservative or in between. The mathematical formulation of $\omega_t$ is as:
\begin{equation}\label{eq:state_def}
\omega_t = \|I'(\textbf{x}^{s,p}_t-\textbf{x}^{e,p}_{t})\|_2,
\end{equation}
where $\textbf{x}^{s,p}_t=[\begin{array}{cccc}x^{s,p}_t & y^{s,p}_t & \theta^{s,p}_t & v^{s,p}_t\end{array}]^T$ denotes the surrounding vehicle's predicted driving state at time step $t$ and $\textbf{x}^{e,p}_t$ is the ego vehicle's predicted driving state made by the surrounding vehicle defined similarly with $\textbf{x}^{s,p}_t$. The definition of $I'$ is as
$$I'=\left[\begin{array}{cccc}
\frac{1}{a_s} & 0 & 0 & 0\\
0 & \frac{1}{b_s} & 0 & 0
\end{array}\right],$$
where $a_s$, $b_s$ are semi-major axis and semi-minor axis of the ellipse representing the surrounding vehicle respectively. The definition (\ref{eq:state_def}) shows that $\omega$ is a scalar and the smaller it is, the less conservative (more aggressive) the ego vehicle appears to the surrounding vehicle. For instance, in the extreme scenario where the distance $d_{t}$ is almost $0$ and $\omega$ equals to $0$ which means the ego vehicle totally does not care about the surrounding vehicle so that its behavior becomes significantly influential to the surrounding vehicle. An opposite case happens when $d_{t}$ is approaching infinity and $\omega$ also goes to infinity, the ego vehicle then cannot bring any influence to the surrounding vehicle and the surrounding vehicle thus has no incentive to consider what the ego vehicle will do in the future. Of course, these scenarios are impossible in the practical driving system. However, they effectively demonstrate the essence of $\omega_t$.

In summary, in the Bayesian persuasive driving process, the ego vehicle aims at finding an optimal signal determined by its driving behavior. The general Bayesian persuasion optimization (\ref{eq:bay_per}) can be reformulated as
\begin{eqnarray}\label{eq:bay_per_car}
\min_{\pi(\textbf{x}^e_t),\textbf{u}^e_t} & \int_{\textbf{x}^e_t} \pi(\textbf{x}^e_t)\min_{\textbf{x}^s_t} \int_{\omega_t} \mu(\omega_t|\textbf{x}^e_t)c(\omega_t,\textbf{x}^s_t,\textbf{x}^e_t,\textbf{u}^e_t)
\end{eqnarray}
for the interactive driving scenario with the expectation term expanded, where
$$\textbf{u}^e_t=[\begin{array}{cc}a^e_t & \delta^e_t\end{array}]^T$$
denotes the ego vehicle's control input at time instant $t$ including acceleration $a^e_t$ and steering angle $\delta^e_t$. In order to avoid shortsighted non-optimal behavior, a receding time horizon is introduced:
\begin{eqnarray}\label{eq:bay_per_car_horizon}
\min_{\pi(\textbf{x}^e_{t|t_0}),\textbf{u}^e_{t|t_0}} && \sum_{t=t_0}^{t_0+N}\int_{\textbf{x}^e_{t|t_0}} \pi(\textbf{x}^e_{t|t_0})\times\nonumber\\
&&\min_{\textbf{x}^s_{t|t_0}} \int_{\omega_{t|t_0}} \mu(\omega_{t|t_0}|\textbf{x}^e_{t|t_0})c(\omega_{t|t_0},\textbf{x}^s_{t|t_0},\textbf{x}^e_{t|t_0},\textbf{u}^e_{t|t_0}),
\end{eqnarray}
where $N$ is the optimization horizon length, $t_0$ is the current time instant, $(\bullet)_{t|t_0}$ denotes prediction of variable for $t$ made at $t_0$.

\subsection{Gaussian assumption}
The optimization (\ref{eq:bay_per_car}) is intractable since the decision variable is a probability distribution in continuous space. In order to make (\ref{eq:bay_per_car}) solvable, the Gaussian assumption is made so that each probability distribution can be described by a mean and a variance in an exponential form. Therefore, the decision variable of (\ref{eq:bay_per_car}) is reduced from a complicated probability distribution to a vector-valued mean and a covairance matrix. Moreover, the exponential form of Gaussian distribution also facilitates the next integral approximation step. With the Gaussian assumption, the probability distributions $\pi$ and $\mu$ in the original problem (\ref{eq:bay_per_car_horizon}) can be explicitly written as
\begin{eqnarray}
&&\pi:\textbf{x}_{t|t_0}^e\sim \mathcal{N}(\hat{\textbf{x}}^e_{t|t_0}, \Sigma_{\textbf{x}_{t|t_0}^e}),\nonumber\\
&&\mu:\omega_{t|t_0}|\textbf{x}^e_{t|t_0}\sim \mathcal{N}(\hat{\omega}_{t|t_0}|\textbf{x}^e_{t|t_0}, \Sigma_{\omega_t|\textbf{x}^e_{t|t_0}}),
\end{eqnarray}
where $x\sim\mathcal{N}(\hat{x}, \Sigma)$ means that random variable $x$ has the Gaussian distribution with mean of $\hat{x}$ and covariance of $\Sigma$.

Then the Bayesian persuasion game can be reorganized as:
\begin{eqnarray}\label{eq:bay_per_car_g}
\min_{\hat{\textbf{x}}^e_{t|t_0}, \textbf{u}^e_{t|t_0}}&&\sum_{t=t_0}^{t_0+N}\int_{\textbf{x}^e_{t|t_0}} G_{\textbf{x}^e_{t|t_0}}(\hat{\textbf{x}}^e_{t|t_0}, \Sigma_{\textbf{x}^e_{t|t_0}})\times\nonumber\\
&&\min_{\textbf{x}^s_{t|t_0}} \int_{\omega_{t|t_0}} \left[G_{\omega_{t|t_0}}(\hat{\omega}_{t|t_0}|\textbf{x}^e_{t|t_0}, \Sigma_{\omega_{t|t_0}|\textbf{x}^e_{t|t_0}})\times\right.\nonumber\\
&&\left.\quad\quad\quad\quad c(\omega_{t|t_0},\textbf{x}^s_{t|t_0},\textbf{x}^e_{t|t_0},\textbf{u}^e_{t|t_0})\right]
\end{eqnarray}
where the expectations of ego vehicle driving states $\hat{\textbf{x}}_{t|t_0}^e (t=t_0,\cdots,t_0+N)$ and control inputs $\textbf{u}_{t|t_0}^e (t=t_0,\cdots,t_0+N)$ are new decision variables,
$$
G_x(\hat{x},\Sigma)=\frac{exp(-\frac{1}{2}(x-\hat{x})^T\Sigma^{-1}(x-\hat{x}))}{\sqrt{(2\pi)^k|\Sigma|}},
$$
denotes the density function of Gaussian distribution $\mathcal{N}(\hat{x},\Sigma)$ and $k$ is the dimension of $x$. $\hat{\omega}_{t|t_0}|\textbf{x}^e_{t|t_0}$ follows the same definition of $\omega$ as in (\ref{eq:state_def}):
\begin{equation}\label{eq:omega_hat}
\hat{\omega}_{t|t_0}|\textbf{x}^e_{t|t_0}=\|I'(\textbf{x}^{s,p}_t-\textbf{x}^e_{t|t_0})\|_2,
\end{equation}
where the ego vehicle's driving state prediction $\textbf{x}^{e,p}_t$ is replaced by $\textbf{x}^{e}_{t|t_0}$.
According to the definition (\ref{eq:omega_hat}), the surrounding vehicle's expected impression of the ego vehicle given his driving behavior is dependent on how the ego vehicle will influence its original driving plan.

\subsection{Cost function}
The cost function $c$ in (\ref{eq:bay_per_car_g}) is defined in exponential form as
\begin{eqnarray}\label{eq:cost}
c(\omega_{t|t_0},\textbf{x}^s_{t|t_0},\textbf{x}^e_{t|t_0},\textbf{u}^e_{t|t_0})&=&exp\left((1+\omega_{t|t_0})\|\textbf{x}^s_{t|t_0}-\textbf{x}^{s,p}_{t|t_0}\|_{W_1}^2\right.\nonumber\\
&&+\|\textbf{x}^e_{t|t_0}-\textbf{x}^e_g\|^2_{W_2}+\|\Delta\textbf{u}^{e,T}_{t|t_0}\|_{W_3}\nonumber\\
&&\left.-\|\textbf{x}^s_{t|t_0}-\textbf{x}^e_{t|t_0}\|^2_{W_4}\right),
\end{eqnarray}
where $\textbf{x}^e_g$ denotes the ego vehicle's desired goal state, $\Delta\textbf{u}^{e}_{t|t_0}=\textbf{u}^{e}_{t|t_0}-\textbf{u}^{e}_{t-1|t_0}$ represents the change of ego vehicle's control input and $W_1$, $W_2$, $W_3$ and $W_4$ are positive definite penalty matrices. It is required that $W_1-W_4\succeq0$ in order to guarantee the existence of a minimum for the cost function $c$ with regard to $\textbf{x}_{t|t_0}^s$. In the rest of the paer, $W_1$ and $W_4$ are set to be equal to $w_1I$ and $w_4I$ respectively, where $I$ denotes the identity matrix, $w_1$ and $w_4$ are scalars.

The first term in the cost function (\ref{eq:cost}) represents the surrounding vehicle's preference of tracking his original driving plan, the second and third terms are driving the ego vehicle to his goal and penalizing the input change in order to achieve comfortable driving experience and the last term represents the surrounding vehicle's aversion of risk, i.e., the preference to keep a certain distance from the ego vehicle. According to the definition (\ref{eq:cost}), when the perceived conservativeness of the ego vehicle $\omega_{t|t_0}$ is lower, the surrounding vehicle will be inclined to focus more on the safety instead of sticking to his original plan. Otherwise, when $\omega_{t|t_0}$ is higher, meaning that the surrounding vehicle is more likely to treat the ego vehicle as a conservative agent, it will be intuitive for him to pursue a more selfish behavior.

The cost function $c$ in (\ref{eq:cost}) shows that the surrounding vehicle's action is dependent on two factors, i.e., the world state of the Bayesian game $\omega_{t|t_0}$ and the penalty matrix $W_1$. The penalty matrix $W_1$ represents the interacting vehicle's driving characteristics, which can only be recognized but not controlled.

Note that although the cost function (\ref{eq:cost}) is intended for two vehicle interaction scenario, the framework can be extended to multiple vehicle interaction case by including more surrounding vehicles in the cost function definition.

\subsection{Constraints}
In the interactive driving application, besides the Bayesian cost function in (\ref{eq:bay_per_car_g}), certain constraints need to be handled including model dynamics, control input saturation and safety constraint.

\subsubsection{Vehicle dynamics}
In this paper, the bicycle model \cite{bicycle_model} is adopted to describe the vehicle dynamics as follows:
\begin{eqnarray}
&&x^e_{t+1|t_0} {=} x^e_{t|t_0} {+} T_s v^e_{t|t_0} \cos\left(\theta^e_{t|t_0} {+} \tan^{-1}(\frac{L_r}{L}\tan\delta^e_{t|t_0})\right)\nonumber\\
&&y^e_{t+1|t_0} {=} y^e_{t|t_0} {+} T_s v^e_{t|t_0} \sin\left(\theta^e_{t|t_0} {+} \tan^{-1}(\frac{L_r}{L}\tan\delta^e_{t|t_0})\right)\nonumber\\
&&\theta^e_{t+1|t_0} {=} \theta^e_{t|t_0} {+}T_s v^e_{t|t_0} \frac{\tan\delta^e_t}{L} \cos\left(\tan^{-1}(\frac{L_r}{L}\tan\delta^e_{t|t_0})\right)\nonumber\\
&&v^e_{t+1|t_0} {=} v^e_{t|t_0} {+}T_s a^e_{t|t_0},\label{eq:state}
\end{eqnarray}
where $T_s$ is the sampling time, $t$ is the time index, $L_r$, $L_f$ and $L{=}L_r{+}L_f$ denote the dimension parameters of the vehicle which are the vehicle's rear, front and full length respectively.

The model equations (\ref{eq:state}) can be summarized as
\begin{equation}
\textbf{x}^e_{t+1|t_0} = f(\textbf{x}^e_{t|t_0}, \textbf{u}^e_{t|t_0}).
\end{equation}

\subsubsection{Safety constraint}
Another critical constraint for autonomous driving is the guarantee of safety. For the static obstacles including parking vehicles and lane boundaries, the following constraint is defined:
\begin{equation}
y_{min}\le y^{e,i}_{t|t_0} \le y_{max}, i\in \{1,2,3,4\},
\end{equation}
where $y^{e,i}_{t|t_0}$ denotes the  y coordinate of the vehicle's $i$-th corner at time step $t$, $y_{min}$ and $y_{max}$ represent the lateral position's lower and upper limit respectively. With regard to the moving surrounding vehicle obstacles, the safety constraint is defined as:
\begin{equation}\label{eq:cons_safe}
\|I'(\textbf{x}^{s,p}_t-\textbf{x}^e_{t|t_0})\|_2\ge 1,
\end{equation}
where the surrounding vehicle is described by the same ellipse as in (\ref{eq:state_def}). In addition, another constraint is imposed on the vehicle speed which is
\begin{equation}\label{eq:cons_v}
v_{min}\le v^{e}_{t|t_0} \le v_{max},
\end{equation}
where $v_{min}$ and $v_{max}$ are minimum and maximum speed respectively.

\subsubsection{Saturation constraint}
In addition to the vehicle modeling and safety constraint, the system also needs to be consistent with the physical control saturation constraint described by
\begin{equation}
\underline{\textbf{u}}^e\le \textbf{u}^e_{t|t_0}\le \overline{\textbf{u}}^e,
\end{equation}
where $\underline{\textbf{u}}^e$, $\overline{\textbf{u}}^e$ represent the lower and upper control saturation bound respectively.

\subsection{Integral calculation}
Currently, the Bayesian persuasion is mainly applied in the economy community. The main factor limiting its popularity in other areas is that calculation of expectation for continuous distribution requires computation of integrals as illustrated by (\ref{eq:bay_per_car_g}). Although Gaussian assumption grants the reduction of decision variable from a distribution to a vector, calculation of integrals is still challenging, especially for the autonomous driving problem with high dimension involved. In this subsection, through introducing several approximations, the original cost fucntion (\ref{eq:bay_per_car_horizon}) is reformulated into a tractable form. By an abuse of notation,  the subscript $t|t_0$ is replaced by $t$ in this subsection.

First consider the integral with regard to $\omega_{t}$:
\begin{equation}\label{eq:int_omega}
\int_{\omega_{t}} G_{\omega_{t}}(\hat{\omega}_{t}|\textbf{x}^e_{t}, \Sigma_{\omega_{t}|\textbf{x}^e_{t}})c(\omega_{t},\textbf{x}^s_t,\textbf{x}^e_{t},\textbf{u}^e_{t}).
\end{equation}
As $\omega_{t}$ is a scalar variable, we can directly calculate the result of (\ref{eq:int_omega}) as
\begin{eqnarray}\label{eq:sub_min}
&&\int_{\omega_t} G_{\omega_t}(\hat{\omega}_t|\textbf{x}^e_t, \Sigma_{\omega_t|\textbf{x}^e_t})c(\omega_t,\textbf{x}^s_t,\textbf{x}^e_t,\textbf{u}^e_t)\nonumber\\
&=& \epsilon exp\left(\frac{1}{2}\Sigma_{\omega_t|\textbf{x}^e_t}[(\textbf{x}^s_t-\textbf{x}^{s,p}_t)^TW_1(\textbf{x}^s_t-\textbf{x}^{s,p}_t)]^2\right)\times\nonumber\\
&&c(\hat{\omega_t}|\textbf{x}^e_t,\textbf{x}^s_t,\textbf{x}^e_t,\textbf{u}^e_t),
\end{eqnarray}
where
$$
\epsilon=\Psi\left(\frac{\hat{\omega}_t|\textbf{x}^e_t+\Sigma_{\omega_t|\textbf{x}^e_t}(\textbf{x}^s_t-\textbf{x}^{s,p}_t)^TW_1(\textbf{x}^s_t-\textbf{x}^{s,p}_t)}{\Sigma_{\omega_t|\textbf{x}^e_t}^{\frac{1}{2}}}\right),
$$
$$
\Psi(x)=\frac{1}{2}(1+erf(x/\sqrt{2}))
$$
is the cumulative distribution function of a standard Gaussian distribution and $erf(\bullet)$ is the error function. According to the safety constraint (\ref{eq:cons_safe}), $\hat{\omega}_t|x^e_t\ge 1$ always holds. Hence when the variance $\Sigma_{\omega_t|\textbf{x}^e_t}$ is chosen to be less than $0.2$, we have
\begin{equation}\label{eq:approx_epsilon}
1\ge\epsilon\ge \Psi(\sqrt{5})=\frac{1}{2}(1+erf(\sqrt{2.5}))= 0.9873,
\end{equation}
resulting in $\epsilon\approx 1$ due to $erf(\bullet)$'s S-shape property.

Then, the solution to the sub-problem
$$
\min_{\textbf{x}^s_t} \int_{\omega_t} G_{\omega_t}(\hat{\omega}_t|\textbf{x}^e_t, \Sigma_{\omega_t|\textbf{x}^e_t})c(\omega_t,\textbf{x}^s_t,\textbf{x}^e_t,\textbf{u}^e_t)
$$
can be obtained as
\begin{eqnarray}\label{eq:xh_sol}
\textbf{x}^{s*}_t&\approx& \left((1+\Sigma_{\omega_t|\textbf{x}^e_t}\gamma_t+\hat{\omega}|\textbf{x}^e_t)W_1-W_4\right)^{-1}*\nonumber\\
&& \left((1+\Sigma_{\omega_t|\textbf{x}^e_t}\gamma_t+\hat{\omega}|\textbf{x}^e_t)W_1\textbf{x}^{s,p}_t-W_4\textbf{x}^e_t\right)
\end{eqnarray}
via setting the derivative equal to zero, where
$$
\gamma_t=(\textbf{x}^s_{t-1}-\textbf{x}^{s,p}_{t-1})^TW_1(\textbf{x}^s_{t-1}-\textbf{x}^{s,p}_{t-1}),
$$
$\textbf{x}^{s,p}_{t-1}$ is the previous prediction and $\textbf{x}^s_{t-1}$ denotes the surrounding vehicle's previous state. 

Substituting (\ref{eq:xh_sol}) into (\ref{eq:sub_min}) obtains
\begin{eqnarray}\label{eq:sol_int_1}
&&\min_{\textbf{x}^s_t}\int_{\omega_t} G_{\omega_t}(\hat{\omega}_t|\textbf{x}^e_t, \Sigma_{\omega_t|\textbf{x}^e_t})c(\omega_t,\textbf{x}^s_t,\textbf{x}^e_t)\nonumber\\
&\approx&exp\{(\textbf{x}^e_t-\textbf{x}^e_g)^TW_2(\textbf{x}^e_t-\textbf{x}^e_g)+\textbf{u}^{e,T}_tW_3\textbf{u}^e_t+k_1\hat{\omega}_t|\textbf{x}^e_t+k_2\nonumber\\
&&+\frac{k_3}{[w_1\hat{\omega}_t|\textbf{x}^e_t+((1+\Sigma_{\omega_t|\textbf{x}^e_t}\gamma_t)w_1-w_4)]}\nonumber\\
&&+\frac{k_4}{[w_1\hat{\omega}_t|\textbf{x}^e_t+((1+\Sigma_{\omega_t|\textbf{x}^e_t}\gamma_t)w_1-w_4)]^2}\},
\end{eqnarray}
where $k_1$, $k_2$, $k_3$ and $k_4$ are constant scalars dependent on $w_1$, $w_4$, $\Sigma_{\omega_t}|\textbf{x}^e_t$ and $\gamma_t$. Furthermore, it can be guaranteed that $k_3$ is always negative and $k_4$ is always positive so that there exists a upper bound $\overline{k}$ for the last two terms in (\ref{eq:sol_int_1}) as $\hat{\omega}_t|\textbf{x}^e_t\ge0$. 

Thus the original optimization (\ref{eq:bay_per_car_g}) is changed to
\begin{eqnarray}
\min_{\hat{\textbf{x}}^e_t,\textbf{u}^e_t} \int_{\textbf{x}^e_t} G_{\textbf{x}^e_t}(\hat{\textbf{x}}^e_t,\Sigma_{\textbf{x}^e_t})exp\{(\textbf{x}^e_t-\textbf{x}^e_g)^TW_2(\textbf{x}^e_t-\textbf{x}^e_g)\nonumber
\end{eqnarray}
\begin{equation}\label{eq:int_2}
+\textbf{u}^{e,T}_tW_3\textbf{u}^e_t+k_1\hat{\omega}_t|\textbf{x}^e_t\},\quad\quad\quad\quad\quad\quad
\end{equation}
where the constant terms $k_2$ and $\overline{k}$ are omitted. 

The integral term in (\ref{eq:int_2}) can be compactly written as
$$
\int_{\textbf{x}^e_t} exp(C(\textbf{x}^e_t)),
$$
where $C(\textbf{x}^e_t)$ is in quadratic form. Since the variable $\textbf{x}^e_t$ is a vector, a mathematical approximation technique based on Lapace's method is utilized to eliminate the high-dimensional integral \cite{dragan2012formalizing}. Via taking second order Taylor seires expansion around $\textbf{x}^{e*}_t{:=}arg\min_{\textbf{x}^e_t}C(\textbf{x}^e_t)$ where $\nabla C(\textbf{x}^{e*}_t)$ equals to $0$, the following equation can be obtained:
\begin{eqnarray}
C(\textbf{x}^e_t)&\approx& C(\textbf{x}^{e*}_t)+\nabla C(\textbf{x}^{e*}_t)^T(\textbf{x}^e_t-\textbf{x}^{e*}_t)\nonumber\\
&&+\frac{1}{2}(\textbf{x}^e_t-\textbf{x}^{e*}_t)^T\nabla^2C(\textbf{x}^{e*}_t)(\textbf{x}^e_t-\textbf{x}^{e*}_t)\nonumber\\
&=& C(\textbf{x}^{e*}_t)+\frac{1}{2}(\textbf{x}^e_t-\textbf{x}^{e*}_t)^T\nabla^2C(\textbf{x}^{e*}_t)(\textbf{x}^e_t-\textbf{x}^{e*}_t).
\end{eqnarray}
Then the integral term $\int_{\textbf{x}^{e}_t}exp(C(\textbf{x}^e_t))$ can be approximated as
\begin{eqnarray}
\int_{\textbf{x}^e_t}exp(C(\textbf{x}^e_t))&\approx& k_e exp(C(\textbf{x}^{e*}_t)),
\end{eqnarray}
where $k_e$ is a constant determined by the Hessian $\nabla^2C(\textbf{x}^{e*}_t)$. Then the cost function in (\ref{eq:int_2}) can be reformulated as
\begin{eqnarray}\label{eq:sol_int_2}
\min_{\hat{\textbf{x}}^e_t, \textbf{u}^e_t} &J_t&= -\frac{1}{2}(\textbf{x}^{e*}_t-\hat{\textbf{x}}^e_t)^T\Sigma_{\textbf{x}^e_t}^{-1}(\textbf{x}^{e*}_t-\hat{\textbf{x}}^e_t)+\textbf{u}^{e,T}_tW_3\textbf{u}^e_t\nonumber\\
&&\quad+(\textbf{x}^{e*}_t-\textbf{x}^e_g)^TW_2(\textbf{x}^{e*}_t-\textbf{x}^e_g)+k_1\hat{\omega}_t|\textbf{x}^{e*}_t.
\end{eqnarray}
Note that $\textbf{x}^{e*}_t$ in (\ref{eq:sol_int_2}) is a linear combination of $\hat{\textbf{x}}^e_t$, $\textbf{x}^e_g$ and $\textbf{x}^{s,p}_t$, thus the cost function in (\ref{eq:sol_int_2}) is a quadratic function with regard to $\hat{\textbf{x}}^e_t$.

\textbf{Remark 1}: Note that the accurate distribution of $\omega_t|\textbf{x}^e_t$ is supposed to be truncated Gaussian instead of standard Gaussian as $\omega_t$ is always positive due to defintion in (\ref{eq:state_def}). The strict formula of $\omega_t|\textbf{x}^e_t$'s density function is
\begin{equation}\label{eq:tr_g}
f(\omega|\textbf{x}^e_t)=\frac{G_{\omega_t}(\hat{\omega_t}|\textbf{x}^e_t, \Sigma_{\omega_t})}{\Sigma_{\omega_t|\textbf{x}^e_t}^{\frac{1}{2}}\left(1-\Psi\left(-\frac{\hat{\omega_t}|\textbf{x}^e_t}{\Sigma_{\omega_t|\textbf{x}^e_t}^{\frac{1}{2}}}\right)\right)}.
\end{equation}
According to the safety consrtaint (\ref{eq:cons_safe}) and utilizing the same approximation as for $\epsilon$ in (\ref{eq:approx_epsilon}), the denominator of (\ref{eq:tr_g}) can be treated as a constant and hence the usage of Gaussian distribution in (\ref{eq:bay_per_car_g}) is validated.

\subsection{Bayesian persuasive optimization}
With the integral calculation result from the previous subsection, the solvable form of the original Bayesian persuasive driving optimization (\ref{eq:bay_per_car_horizon}) can be summarized as
\begin{eqnarray}\label{eq:bay_op}
\min_{\hat{\textbf{x}}^e_{t|t_0},\textbf{u}^e_{t|t_0}}& \sum_{t=t_0}^{t_0+N-1} J_t\\
s.t.& \hat{\textbf{x}}^e_{t+1|t_0} = f(\hat{\textbf{x}}^e_{t|t_0},\textbf{u}^e_{t|t_0}),\nonumber\\
& \|I'(\hat{\textbf{x}}^e_{t|t_0}-\textbf{x}^{s,p}_{t})\|_2\ge 1,\nonumber\\
& \underline{\textbf{u}}^e\le \textbf{u}^e_{t|t_0}\le \overline{\textbf{u}}^e,\nonumber\\
& \hat{\textbf{x}}^e_{t_0|t_0}=\textbf{x}^e_{t_0},v_{min}\le \hat{v}^{e}_{t|t_0} \le v_{max},\nonumber\\
& y_{min}\le \hat{y}^{e,i}_{t|t_0} \le y_{max}, i\in \{1,2,3,4\},\nonumber\\
& t=t_0, \cdots, t_0+N-1,
\end{eqnarray}
where $\textbf{x}^e_{t_0}$ is the ego vehicle's current driving state.

\textbf{Remark 2}: After obtaining the solution to the optimization (\ref{eq:bay_op}) and observing the surrounding vehicle's behavior, more information about the interacting vehicle's characteristics will be revealed so that the covariance and penalty matrices in the cost function should be updated accordingly. With regard to $\Sigma_{\textbf{x}^e_t}$, it is updated based on the surrounding vehicle's confident level of the ego vehicle's behavior, which can be approximated by
\begin{equation}
\alpha=\sum_{t=t_0}^{t_0+N-1}w_t^{\alpha}\|\textbf{x}^e_{t|t_0}-\textbf{x}^e_{t|t_0-1}\|_2,
\end{equation}
where $w_t^{\alpha}$'s are weighting factors satisfying $\sum_{t}w_t^{\alpha}=1$.
According to the definition, smaller $\alpha$ indicates that the ego vehicle's behavior is more consistent with its previous driving plan and thus the surrounding vehicle will be more confident about the information he extract from the ego vehicle's behavior. Therefore, the covariance matrix of probability distribution $\Sigma_{\textbf{x}^r_t}$ should be proportional to $\alpha$.

As shown in the cost function (\ref{eq:cost}), the penalty matrix $W_1$ which reflects the surrounding vehicle's driving characteristics is another influential factor of the algorithm. Similar with the update strategy for $\Sigma_{\textbf{x}^e_t}$, the update of $W_1$ is based on the surrounding vehicle's confidence in his own driving plan, which can be inferred from the following parameter:
\begin{equation}
\beta=\sum_{t=t_0}^{t_0+N-1}w_t^{\beta}\|\textbf{x}^{s,p}_{t|t_0}-\textbf{x}^{s,p}_{t|t_0-1}\|_2,
\end{equation}
where $w_t^{\beta}$'s are weighting factors for $\beta$ satisfying $\sum_{t}w_t^{\beta}=1$. Intuitively, a smaller $\beta$ represents a smaller change of the surrounding vehicle's driving plan, reflecting his more self-centric driving characteristics resulting in a larger $W_1$.

\section{Simulation Results}\label{sec:simu}
In this section, the proposed algorithm's performance is demonstrated by simulations implemented for several driving scenarios including lane changing, lane keeping and intersection crossing. In the simulation, the surrounding vehicles are assumed to follow a model predictive control (MPC) strategy with different safety weights reflecting various driving characteristics, i.e., more cautious/conservative drivers are associated with higher safety weights and more aggressive drivers are represented by lower safety weights. Another assumption is about perfect prediction, which means that the ego vehicle knows exactly what the surrounding vehicle's driving plan is. Although this is a very strong assumption, the algorithm is still necessary as its point is to persuade the surrounding vehicle to change its original plan instead of just identifying it.

\subsection{Simulation setting}
The simulation environment utililzed in this paper is a 1/10 scaled version of the real world and all the vehicles are assumed to be of the same size. The dimension parameters are shown in Table \ref{table:vdp}, where $W_V$ and $W_L$ are width of the vehicle and the lane respectively. The other optimization configurations are shown in Table \ref{table:opc}, where $\underline{u}^e_{\delta}$, $\overline{u}^e_{\delta}$ are lower and upper bound for steering angle and $\underline{u}^e_{a}$, $\overline{u}^e_{a}$ are lower and upper bound for acceleration respectively. 
\begin{table}[!htbp]
	\caption{Dimension parameters}
	\label{table:vdp}
	\begin{center}
		\begin{tabular}{|c|c|c|c|}
			\hline
			$L_f$ & $L_r$ & $W_V$ & $W_L$ \\
			\hline
			0.21(m) & 0.19(m) & 0.19(m) & 0.37(m) \\
			\hline
		\end{tabular}
	\end{center}
\end{table}

\begin{table}[!htbp]
	\caption{Optimization configurations}
	\label{table:opc}
	\begin{center}
		\begin{tabular}{|c|c|c|c|}
			\hline
			$T_s$ & $N$ & $\underline{u}^e_a$ & $\underline{u}^e_{\delta}$ \\
			\hline
			0.1(s) & 30 & -1($m/s^2$) & $-\pi/3$(rad) \\
			\hline
			$\overline{u}^e_a$ & $\overline{u}^e_{\delta}$ & $y_{min}$ & $y_{max}$ \\
			\hline
			1($m/s^2$) & $\pi/3$(rad) & 0(m) & 0.74(m) \\
			\hline
			$v_{min}$ & $v_{max}$ & $a_s$ & $b_s$\\
			\hline
			0(m/s) & 1.5(m/s) & 0.75(m) & 0.35(m)\\
			\hline
		\end{tabular}
	\end{center}
\end{table}

\subsection{Lane changing scenario}
\begin{figure}[thpb]
	\centering
	\includegraphics[width=\columnwidth,height=5cm]{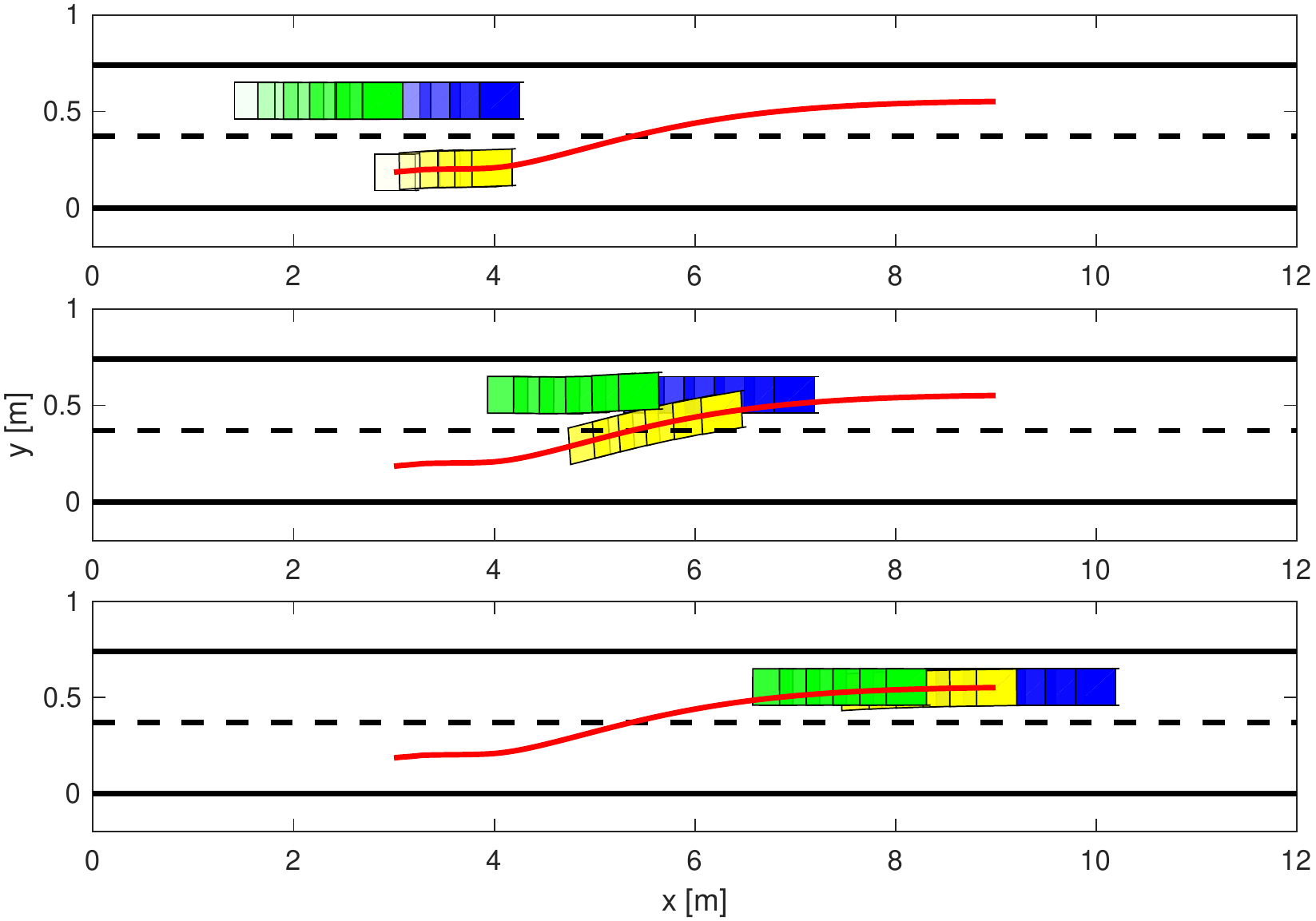}
	\caption{Lane changing scenario with ``nice" surrounding vehicle}
	\label{fig:nice_no_un}
\end{figure}

First consider the scenario where the ego vehicle attempts to change to the neighboring lane with two surrounding vehicles running in it. The front surrounding vehicle is assumed to be always aggressive (safety weight equals to zero) so that it will never yield the ego vehicle and the rear vehicle's safety weight is adjusted to represent different kinds of driver, where large safety weight means nice driver and small or zero safety weight is utilized for iron nerved aggressive driver. The interacting case with ``nice'' driver is shown in Fig.~\ref{fig:nice_no_un}, where the ego vehicle is represented by the yellow rectangle, the front and rear vehicle are plotted as blue and green rectangle respectively. The solid red line represents the ego vehicle's driving  trajectory. In this case, the nice driver decided to decelerate and yield the ego vehicle. The ego vehicle thus took the chance and finished the lane changing task smoothly.

\begin{figure}[thpb]
	\centering
	\includegraphics[width=\columnwidth,height=5cm]{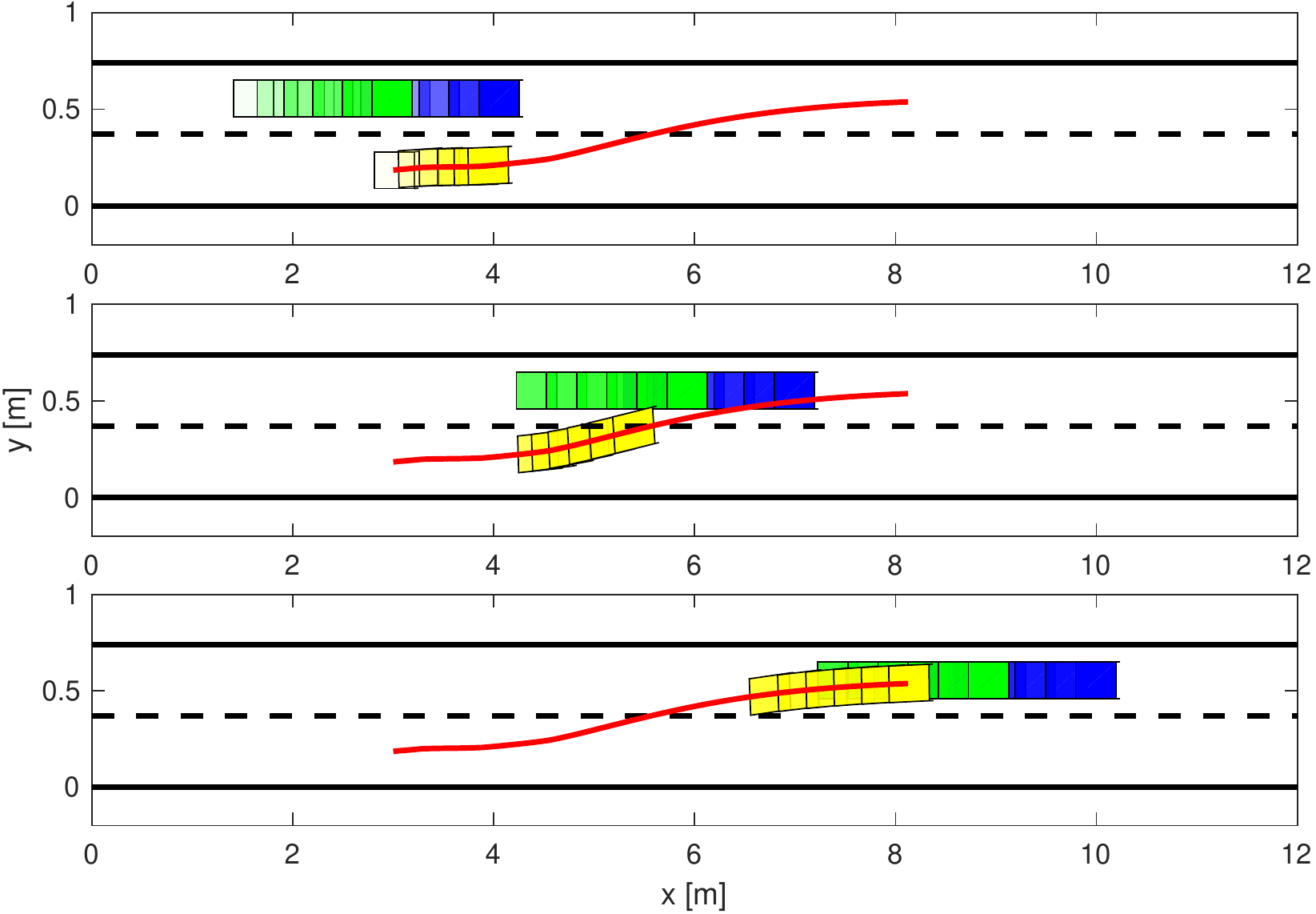}
	\caption{Lane changing scenario with ``tough" surrounding vehicle}
	\label{fig:tough_no_un}
\end{figure}

Figure~\ref{fig:tough_no_un} shows another case with ``tough'' driver who chooses to accelerate and ignore the ego vehicle's lane changing need. In this case, as the rear vehicle refused to yield, the ego vehicle waited until the rear surrounding vehicle passed and completed lane changing later.

\subsection{Lane keeping scenario}
\begin{figure}[thpb]
	\centering
	\includegraphics[width=\columnwidth,height=5cm]{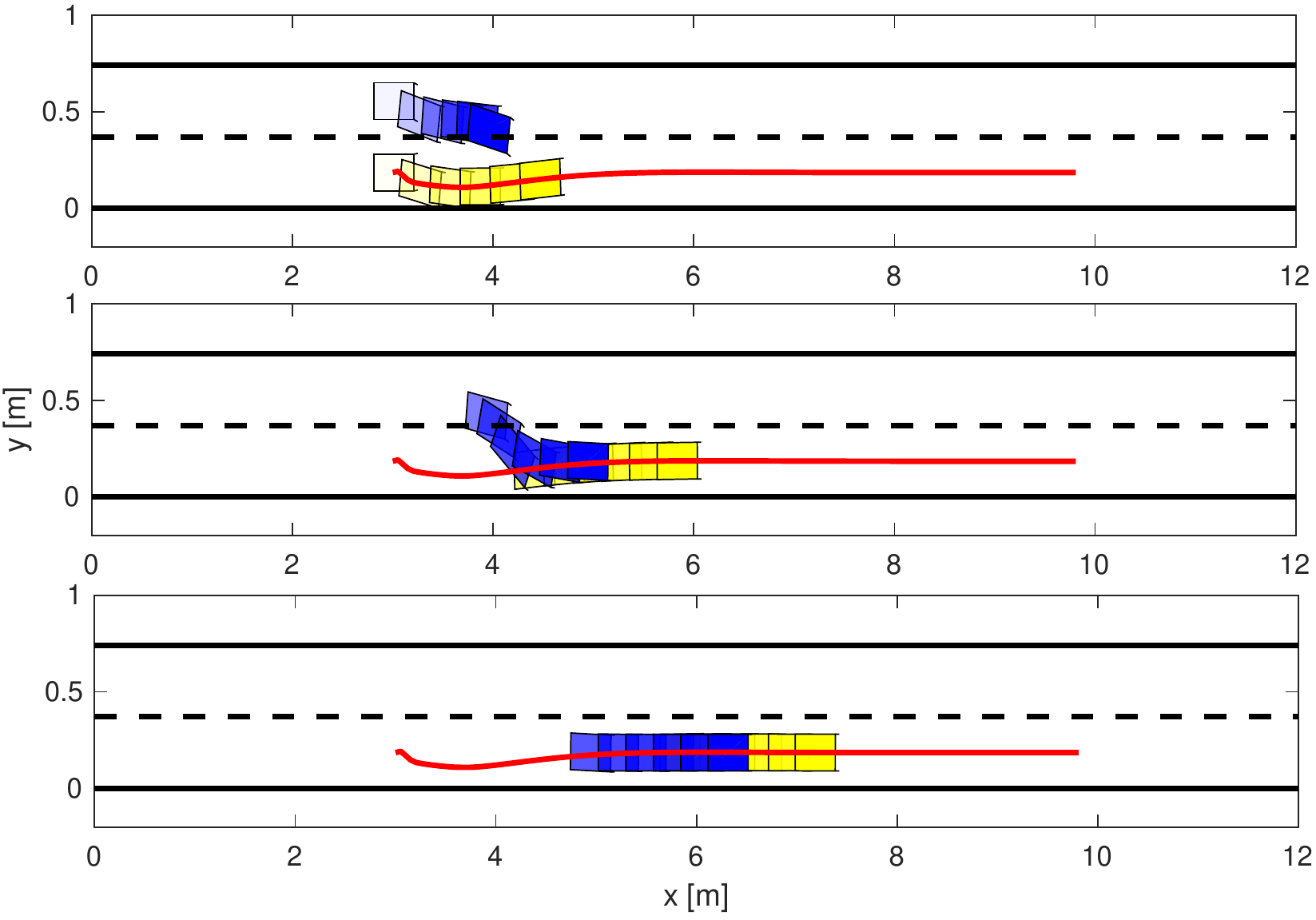}
	\caption{Lane keeping scenario with ``nice" surrounding vehicle}
	\label{fig:lane_keeping_nice}
\end{figure}

\begin{figure}[thpb]
	\centering
	\includegraphics[width=\columnwidth,height=5cm]{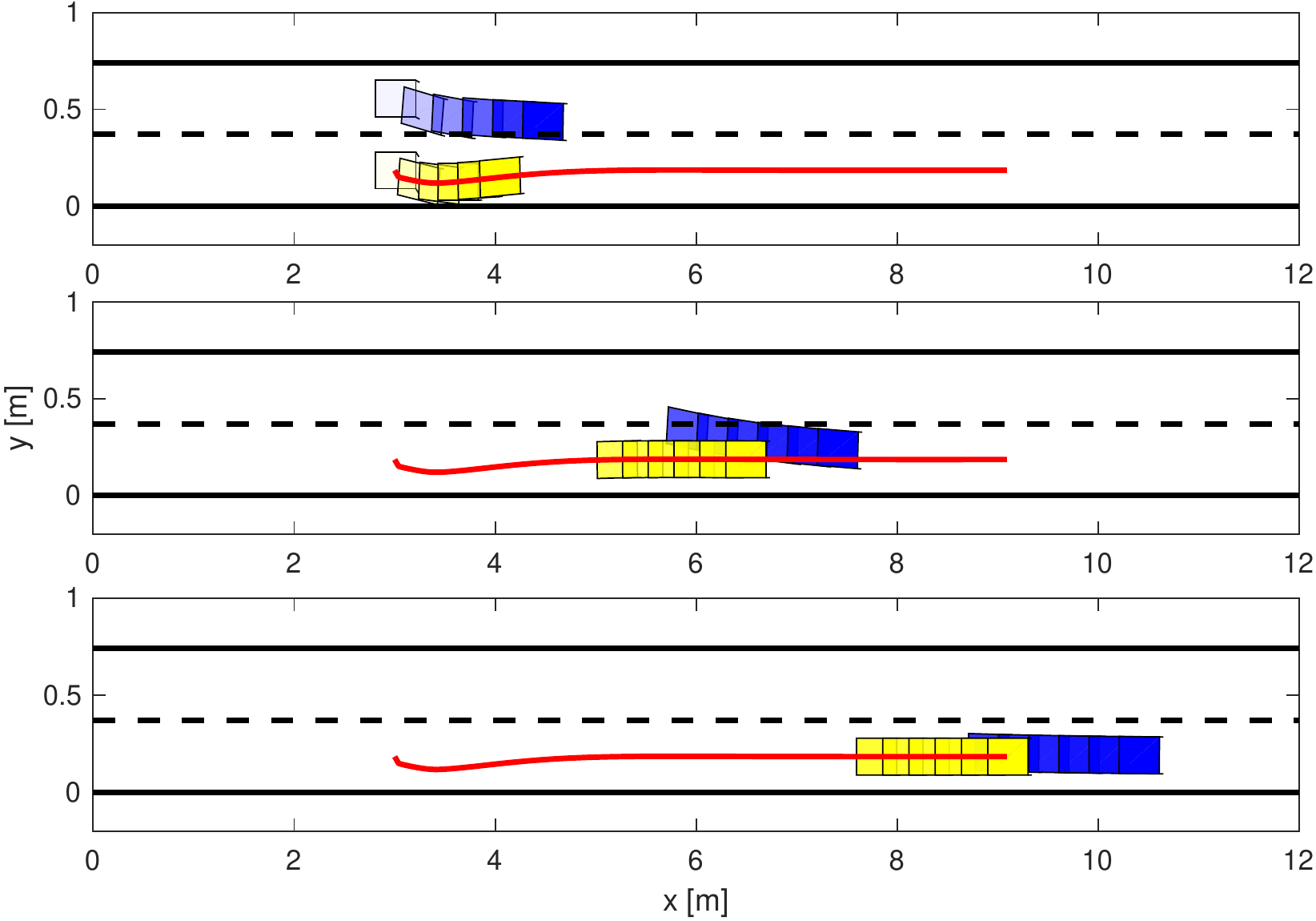}
	\caption{Lane keeping scenario with ``tough" surrounding vehicle}
	\label{fig:lane_keeping_tough}
\end{figure}

\begin{figure}[thpb]
	\centering
	\includegraphics[width=\columnwidth,height=5cm]{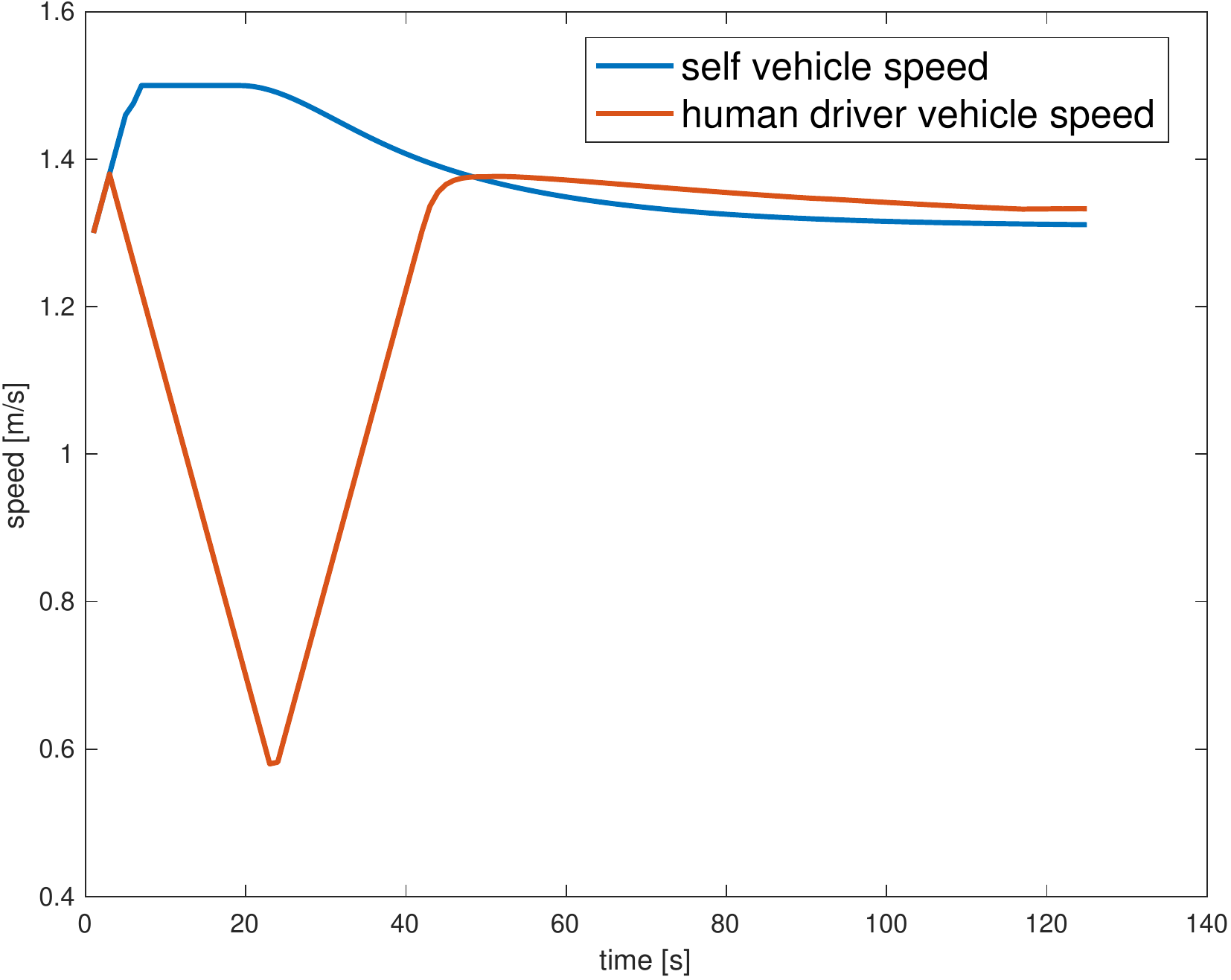}
	\caption{Lane keeping speed profile (``nice" surrounding vehicle case)}
	\label{fig:lane_keep_nice_vel}
\end{figure}

\begin{figure}[thpb]
	\centering
	\includegraphics[width=\columnwidth,height=5cm]{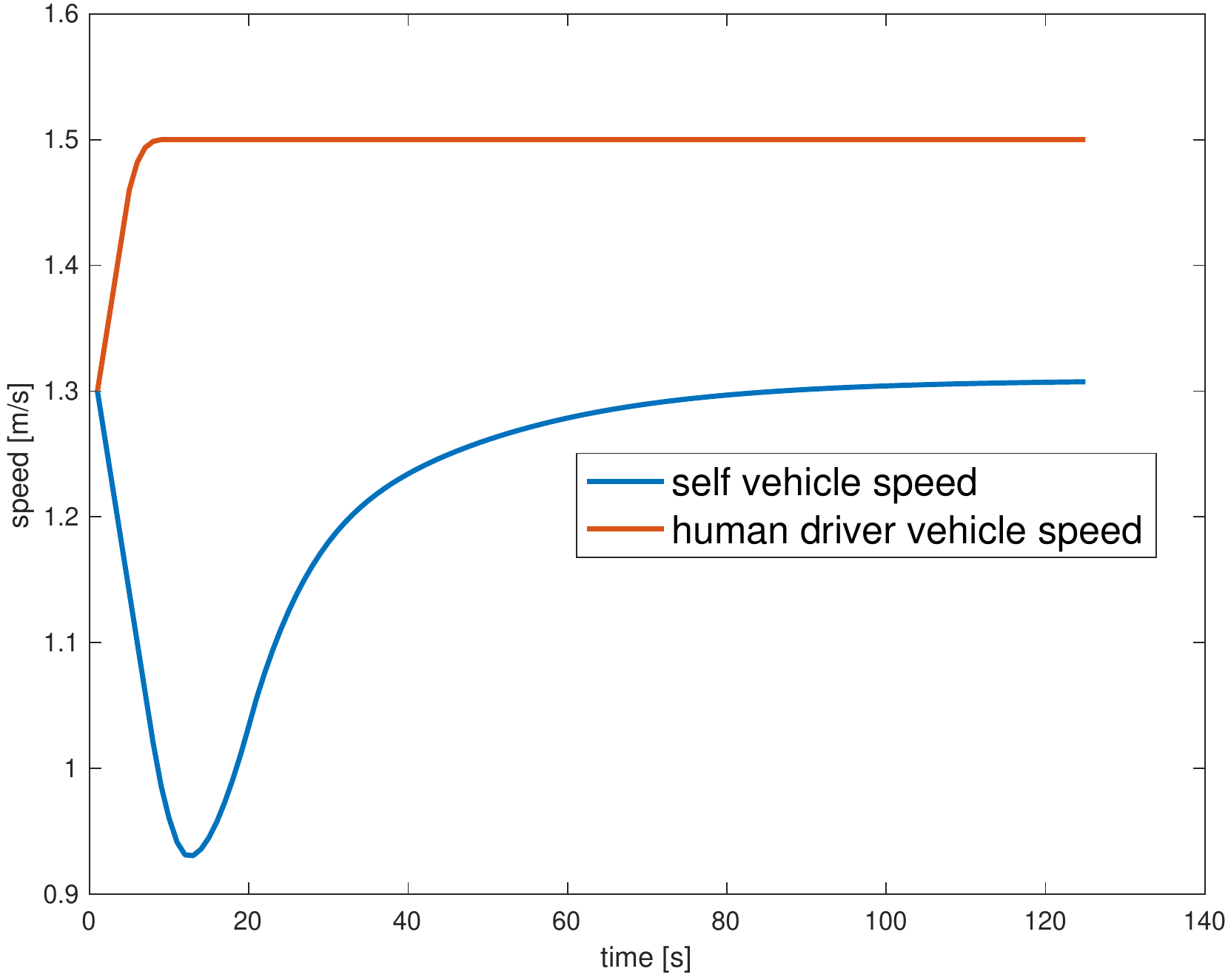}
	\caption{Lane keeping speed profile (``tough" surrounding vehicle case)}
	\label{fig:lane_keep_tough_vel}
\end{figure}

In the lane keeping scenario, the ego vehicle needs to make decision between to yield or not to yield when the surrounding vehicle seeks to merge in. The simulation results Fig.~\ref{fig:lane_keeping_nice} and Fig.~\ref{fig:lane_keeping_tough} illustrate the ego vehicle's intelligent driving behavior when interacting with nice and tough driver respectively. In the case of a nice driver, the ego vehicle decides to ignore its merging request and accelerate to show its intention as shown in Fig.~\ref{fig:lane_keeping_nice}. On the other hand, when interacting with a more aggressive driver, the ego vehicle expresses its intention of yielding by decelerating as shown in Fig.~\ref{fig:lane_keeping_tough}. The speed profile for these two driving situations are shown in Fig.~\ref{fig:lane_keep_nice_vel} and Fig.~\ref{fig:lane_keep_tough_vel} respectively.

\subsection{Intersection crossing scenario}
\begin{figure}[thpb]
	\centering
	\includegraphics[width=\columnwidth,height=3cm]{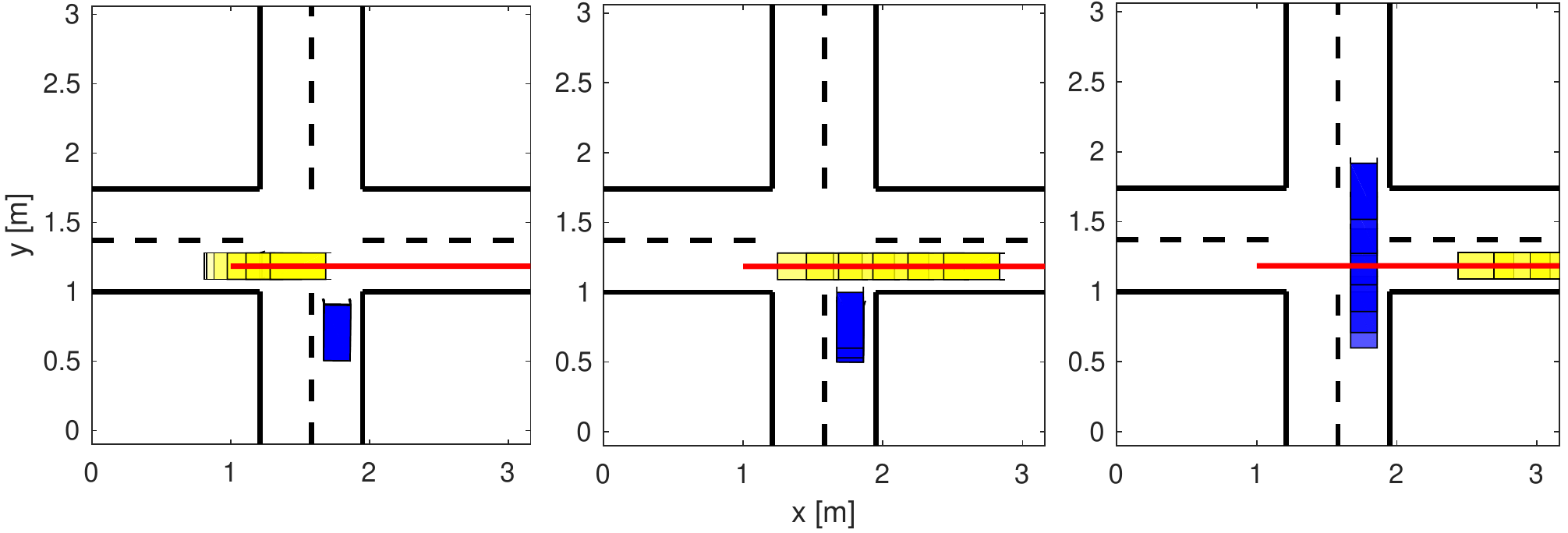}
	\caption{Intersection scenario with ``nice" surrounding vehicle}
	\label{fig:intersection_nice}
\end{figure}

\begin{figure}[thpb]
	\centering
	\includegraphics[width=\columnwidth,height=3cm]{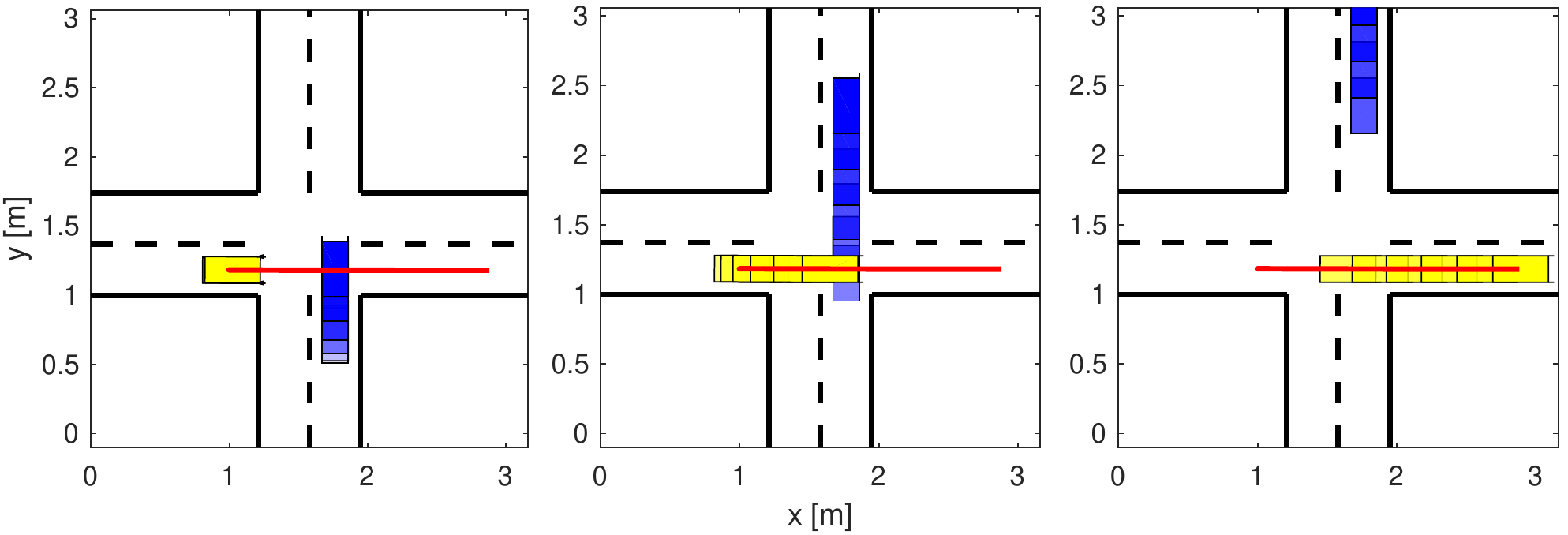}
	\caption{Intersection scenario with ``tough" surrounding vehicle}
	\label{fig:intersection_tough}
\end{figure}

\begin{figure}[thpb]
	\centering
	\includegraphics[width=\columnwidth,height=5cm]{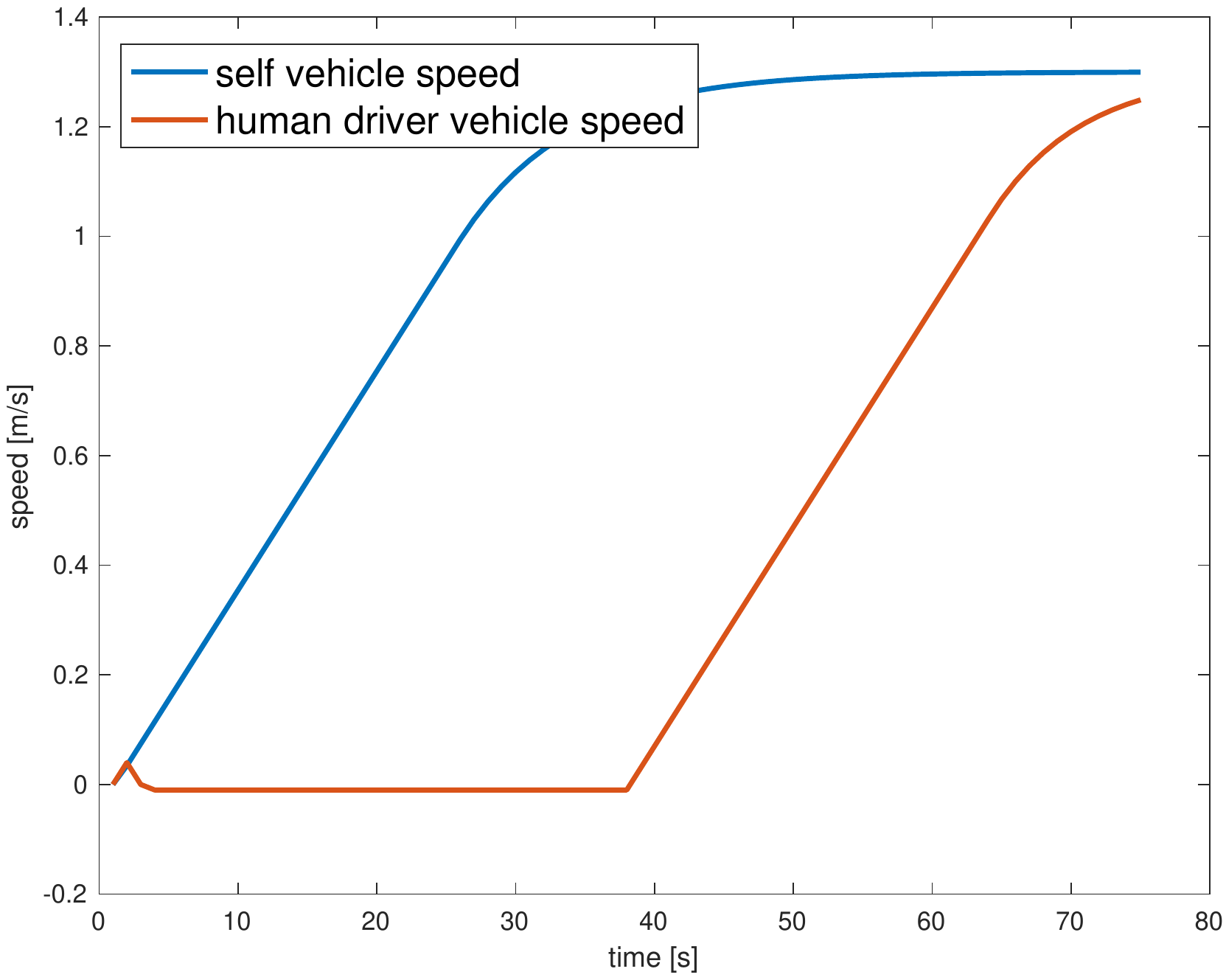}
	\caption{Intersection crossing speed profile (``nice" surrounding vehicle case)}
	\label{fig:intersection_nice_vel}
\end{figure}

\begin{figure}[thpb]
	\centering
	\includegraphics[width=\columnwidth,height=5cm]{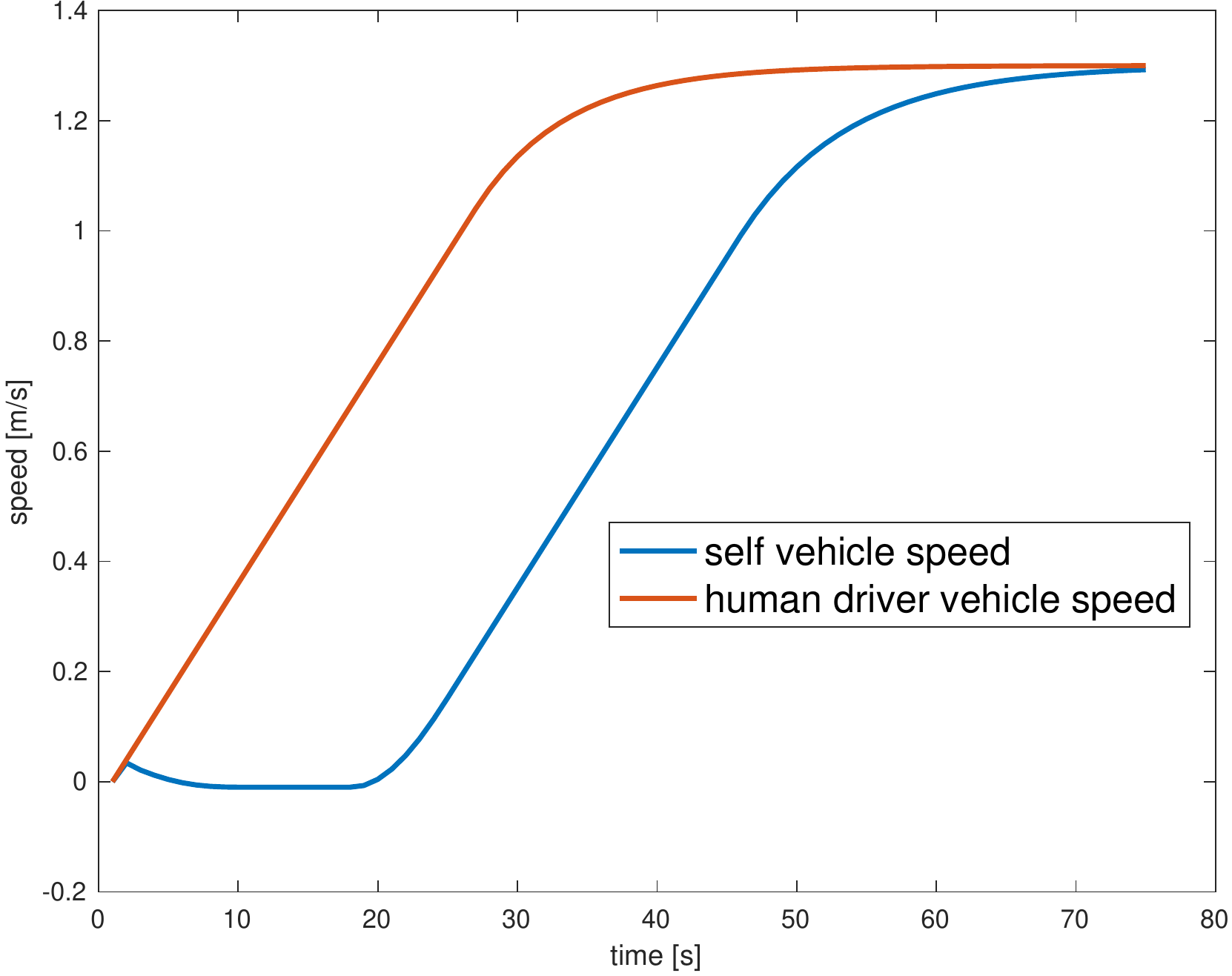}
	\caption{Intersection crossing speed profile (``tough" surrounding vehicle case)}
	\label{fig:intersection_tough_vel}
\end{figure}

Another common driving scenario is intersection crossing, where both the vehicles need to reason about who is supposed to pass first. Fig.~\ref{fig:intersection_nice} and Fig.~\ref{fig:intersection_tough} demonstrate the algorithm's performace when interacting with different kinds of driver. In the figures, the yellow rectangle represent the ego vehicle and the blue rectangle is the surrounding vehicle. It is shown in Fig.~\ref{fig:intersection_nice} (Fig.~\ref{fig:intersection_tough}) that the ego vehicle decides to pass first (second) when meeting a nice (tough) driver, which is consistent with human driving behavior. The vehicles' speed profiles are shown in Fig.~\ref{fig:intersection_nice_vel} and Fig.~\ref{fig:intersection_tough_vel}, which further illustrate the ego vehicle's intention. In the case of nice surrounding vehicle, the ego vehicle keeps accelerating and cross the intersection first. When interacting with a tough surrounding vehicle, the ego vehicle first inches a little and then stops, waiting for the other vehicle to pass.

\section{Conclusion}\label{sec:conc}
In this paper, an optimization-based Bayesian persuasive driving algorithm was proposed. In the persuasion game, the ego vehicle serves as the information sender who attempts to manipulate the surrounding vehicle's (information receiver) posterior belief of the world state in order to achieve a lower cost for both players via providing information about its driving intention. The world state of the Bayesian game was defined to be the surrounding vehicle's impression about the ego vehicle. In the surrounding vehicle's point of view, both the signaling and the belief of the world state are formulated as Gaussian distributions. An integral approximation was applied to reformulate the optimization into a tractable form. As shown by simulation results in several driving scenarios, the ego vehicle is capable of interacting with various types of surrounding vehicles intelligently due to the persuasion signaling strategy.

In our future work, the surrounding vehicle's intention prediction will be studied and incorporated with the proposed Bayesian persuasive algorithm. In addition, a high level decision making controller will also be explored to set appropriate desired goal state for the ego vehicle.

\addtolength{\textheight}{-12cm}   





\bibliographystyle{IEEEtran}
\bibliography{bayesian_persuasion_driving}

\end{document}